\documentclass{notices}

\usepackage{amsmath}
\usepackage{amssymb}
\usepackage{amsthm}
\usepackage{epsfig}
\usepackage{graphicx}
\usepackage{tikz}
\usepackage{tikz-cd}
\usetikzlibrary{arrows, positioning, shapes, decorations.pathreplacing, decorations.markings, calc, matrix}
\usepackage[all]{xy}

\newcommand{\MGC}{\ensuremath{\mathcal{M}_{G_{\mathbb{C}}}}}
\newcommand{\GC}{\ensuremath{G_{\mathbb{C}}}}
\newcommand{\C}{\mathbb{C}}
\newcommand{\R}{\mathbb{R}}

\title{
Higgs bundles --  Recent applications
}

\author{
  Laura P.~Schaposnik  \affil{
    The  author is a professor of mathematics at the University of Illinois at Chicago. Her email address is schapos@uic.edu.
    }
   }

\begin{document}

\maketitle
 
\section*{Introduction}  

This note is dedicated to introducing  Higgs bundles and the Hitchin fibration, with a view towards their appearance  within different branches of mathematics and physics,  focusing in particular on the role played by the \textit{integrable system} structure carried by their moduli spaces. On a compact Riemann surface  $\Sigma$ of genus $g\geq 2$, Higgs bundles  are pairs $(E,\Phi)$ where  \begin{itemize}
 \item $E$   is a holomorphic vector bundle on $\Sigma$;
 \item    the Higgs field $\Phi: E\rightarrow E \otimes K$ is a holomorphic section, for   $K:=T^*\Sigma$.
  \end{itemize}

  Since their origin in the late 80's in work of Hitchin and Simpson, 
Higgs bundles manifest as  fundamental objects which are ubiquitous in contemporary mathematics, and closely related to theoretical physics. For $\GC$ a complex semisimple Lie group, the {\it Dolbeault moduli space} of $G_\mathbb{C}$-Higgs bundles  $\mathcal{M}_{G_\mathbb{C}}$ has a hyperk\"ahler structure, and via different complex structures it  can be seen as different moduli spaces:
 \begin{itemize}
 \item
Via the non-abelian Hodge correspondence developed by Corlette, Donaldson,  Simpson and Hitchin, and in the spirit of Uhlenbeck-Yau's work for compact groups,
 the moduli space   is   diffeomorphic as a real manifold to the {\it De Rahm moduli space} $\mathcal{M}_{dR}$  of flat connections on a smooth complex bundle.
\item Via the Riemann-Hilbert correspondence there is an analytic correspondence 
 between the de Rham   space and the  {\it Betti moduli space}  $\mathcal{M}_{B}$ of
surface group representations.
\end{itemize}
Some prominent examples where these moduli spaces appear in mathematics and physics are:
\begin{itemize}
\item  Through the Hitchin fibration, $\MGC$ give examples of hyperk\"ahler manifolds which are {\it integrable systems},  leading to remarkable applications in physics which we shall discuss later on. 
\item  Building on the work of Hausel and Thaddeus relating Higgs bundles to {\it Langlands duality},     Donagi and Pantev presented $\MGC$ as a fundamental  example of mirror symmetry for CY manifolds.
\item   Within the work of Kapustin and Witten,     Higgs bundles were used to obtain a physical derivation of the {\it geometric Langlands correspondence} through mirror symmetry and soon after, Ng\^{o} found Higgs bundles to be key ingredients when proving the Fundamental Lemma, which led him to the Fields Medal a decade ago. 
 \end{itemize}
 
 Higgs bundles and Hitchin systems have been an increasingly vibrant area, and thus there are several expository articles some of which   we shall refer to: from the  Notices' article  ``What is a Higgs bundle?''  \cite{whatis}, to several graduate notes on Higgs bundles (e.g., the author's recent  \cite{PCMI}), to more advance reviews such as   Ng\^{o}'s 2010 ICM Proceedings article \cite{chau2010endoscopy}. Hoping to avoid repeating material nicely covered in other reviews, whilst still attempting to engage the reader into learning more about the subject, we shall take this opportunity to focus on some of the recent work done by leading young members of the community\footnote{As in other similar reviews, the number of references is limited to twenty, and thus we shall refer the reader mostly to survey articles where precise references can be found.}.

\section*{Higgs bundles}
Higgs bundles arise as solutions to self-dual Yang-Mills equations, a non-abelian generalization of Maxwell's equations  which recurs through  much of modern physics. Recall that instantons are solutions to Yang-Mills self-duality equations in Euclidean 4d space, and when these equations are reduced to  Euclidean 3d space by imposing translational invariance in one dimension, one obtains monopoles as solutions.   
Higgs bundles were introduced by Hitchin in \cite{N1} as solutions of the so-called  {\it Hitchin equations}, the 2-dimensional reduction  of the  Yang-Mills self-duality equations,   given by\begin{eqnarray}
F_A+ [\Phi,\Phi^*]=0,\label{equation1}\\ ~{~}~ ~\overline{\partial}_{A}\Phi=0,\label{equation2}\end{eqnarray}    where $F_A$ is the curvature of     a unitary connection $\nabla_A=\partial_{A}+\overline{\partial}_{A}$ associated to the $\overline{\partial}$-connection   $\overline{\partial}_{A}$ on a principal  $\GC$ bundle $P$. 
Concretely,  principal $\GC$-Higgs bundles  are pairs $(P,\Phi)$, where 
 \begin{itemize}\item $P$ is a principal $\GC$-bundle, and 
 \item $\Phi$  a holomorphic section of ${\rm ad}(P)\otimes K$.
 \end{itemize} 
 Throughout these notes, we shall refer to {\it classical Higgs bundles} as the  Higgs bundles described in the Introduction, and consider $\GC$-Higgs bundles in their vector bundle representation, through which they can be seen as classical Higgs bundles satisfying some extra conditions reflecting the nature of the group $\GC$. For instance when considering $\GC=SL(n,\C)$, a $\GC$-Higgs bundle $(E,\Phi)$ is composed of a holomorphic rank $n$ vector bundle $E$ with trivial determinant $\Lambda^nE\cong \mathcal{O}$, and a Higgs field satisfies ${\rm Tr}(\Phi)=0$, for which we shall write $\Phi \in H^0(\Sigma, {\rm End}_0(E)\otimes K)$.

\smallskip

\noindent \textbf{Example 1}.~Choosing a square root of $K$,    consider the vector bundle $E=K^{1/2}\oplus K^{-1/2}$. Then,   a family of $SL(2,\C)$-Higgs bundles $(E,\Phi_a)$ parametrized by quadratic differentials $a\in H^{0}(\Sigma, K^2)$ is given by 
  \begin{eqnarray}\left(E=K^{1/2}\oplus K^{-1/2}, \Phi_a=\left(\begin{array}{cc}0&a\\1&0\end{array}\right)\right). 
 \label{example} \end{eqnarray}

One may also consider   $G$-Higgs bundles, for $G$ a real form of  $G_\C$, which in turn correspond to the Betti moduli space of representations into $G$. 

\noindent \textbf{Example 2.}   $SL(2,\R)$-Higgs bundles are pairs \[\left(E=L\oplus L^*, \Phi=\left(\begin{array}{cc}0&\beta\\\gamma&0\end{array}\right)\right),\]for $L$ a line bundle on $\Sigma$. Hence, in  {\it Example 1.}  one has a family 
 $(E,\Phi_a)$  of $SL(2,\R)$-Higgs bundles.

\smallskip
   In order to define the moduli space of Higgs bundles, one needs to incorporate the notion of stability. For this, recall that holomorphic vector bundles $E$ on $\Sigma$  are topologically classified by their rank $rk(E)$  and their  degree $deg(E)$, though which one may define their {\it slope} as $\mu(E):= deg(E)/  rk(E).$
Then, a vector bundle $E$ is   {\em stable} ({\em or semi-stable}) if for any proper, non-zero sub-bundle $F\subset E$ one has $\mu(F)<\mu(E)$ (or $\mu(F)\leq \mu(E)$).
 It is {\em  polystable} if it is a direct sum of stable bundles whose slope is  $\mu(E)$.

One can generalize the   stability condition   to Higgs bundles $(E,\Phi)$ by considering $\Phi$-{\em invariant subbundle} $F$ of $E$, a vector subbundle $F$ of $E$ for which $\Phi(F)\subset F\otimes K$. Then,  a Higgs bundle $(E,\Phi)$ is said to be {\em  stable} ({\em semi-stable}) if for each proper $\Phi$-invariant  $F\subset E$ one has $\mu(F)<~\mu(E)~(equiv. \leq)$.   Then, the moduli space $\MGC$ of stable $\GC$-Higgs bundles  up to holomorphic automorphisms of the pairs  can be constructed (also denoted $\mathcal{M}_{Dol}$). Going back to Hitchin's equations, one of the most important characterisations of stable Higgs bundles 
is given in the work of Hitchin and Simpson,  and which carries through to more general settings: 
If a Higgs bundle $(E,\Phi)$ is stable and  $\text{deg} ~E = 0$, then there is a unique unitary connection $A$ on $E$, compatible with the holomorphic structure, satisfying \eqref{equation1}-\eqref{equation2}.

Finally,    Hitchin showed that the moduli space of Higgs bundles  is a hyperk\"ahler manifold with natural symplectic form $\omega$ defined on the infinitesimal deformations
 $(\dot A,\dot \Phi)$ of a Higgs bundle $(E,\Phi)$   by
\begin{equation}\label{ch2:2.1}
 \omega((\dot{A}_{1},\dot{\Phi}_{1}),(\dot{A}_{2},\dot{\Phi}_{2}))=\int_{\Sigma}\text{tr}(\dot{A}_{1}\dot{\Phi}_{2}-\dot{A}_{2}\dot{\Phi}_{1}),
\end{equation}
 where $\dot A \in \Omega^{0,1}(\text{End}_{0} E)$ and $\dot\Phi\in \Omega^{1,0}(\text{End}_0 E)$.
 Moreover, he presented a  natural way of studying the moduli spaces   $\mathcal{M}_{\GC}$ of $ \GC$-Higgs bundles through what is now called {\it the Hitchin fibration}, which we shall consider next.   
  
\newpage
\section*{Integrable systems}

\begin{figure}[h]
\centering
\includegraphics[width=0.4\textwidth]{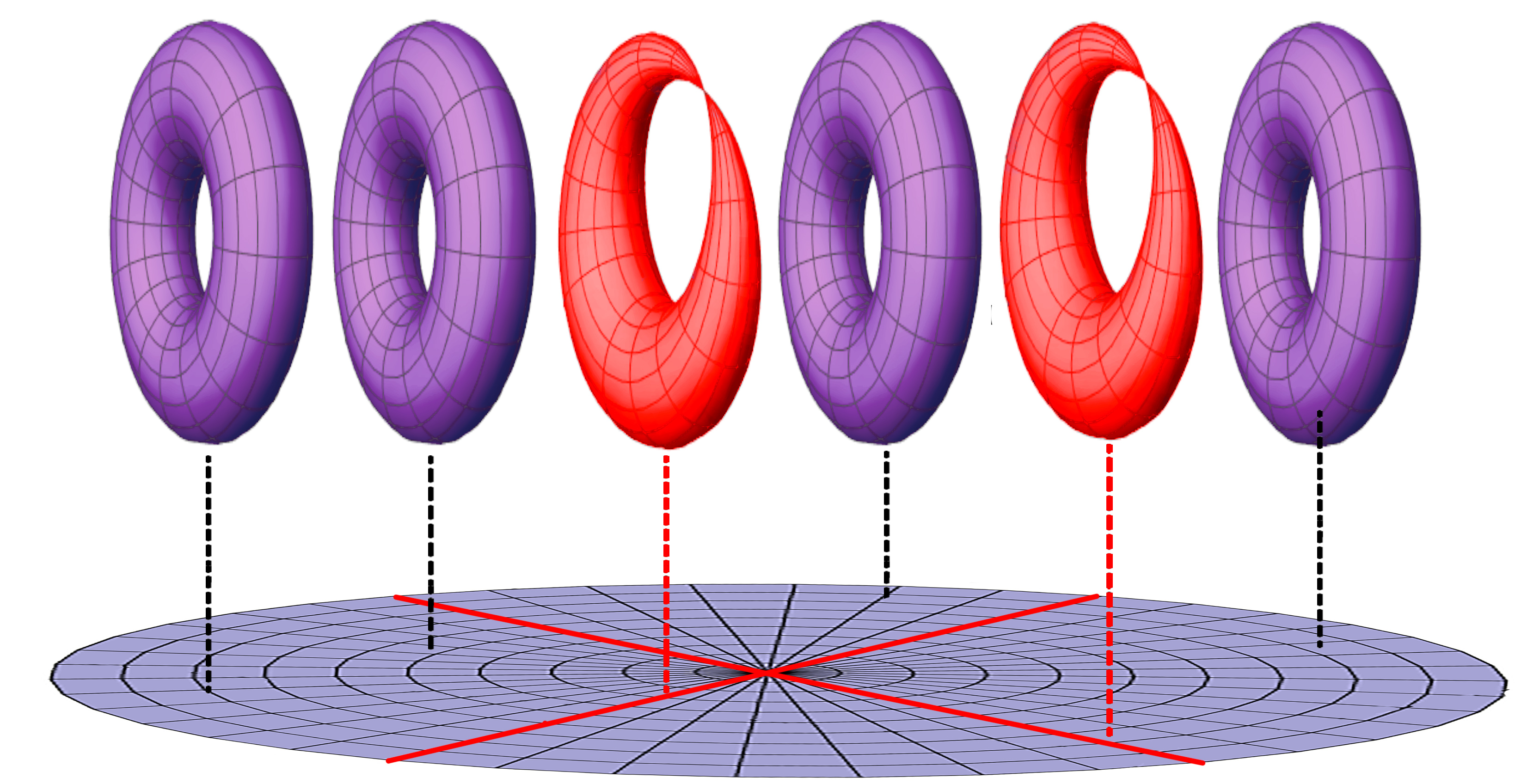}
\caption{The Hitchin fibration.}
\label{fibration} \end{figure}

Given $\{p_{1}, \ldots, p_k\}$ a homogeneous basis for the algebra of invariant polynomials on the Lie algebra  $\mathfrak{g}_{c}$ of $ \GC$,  we denote by $d_{i}$ the degree of $p_i$. The {\it Hitchin fibration}, introduced in   \cite{N2}, is then given by
\begin{eqnarray} h~:~ \mathcal{M}_{ \GC}&\longrightarrow&\mathcal{A}_{ \GC}:=\bigoplus_{i=1}^{k}H^{0}(\Sigma,K^{d_{i}}),  \nonumber  \\
 (E,\Phi)&\mapsto& (p_{1}(\Phi), \ldots, p_{k}(\Phi)).\nonumber
\end{eqnarray} where the map $h$ is referred to as the {\it Hitchin~map}: it is a proper map for any choice of basis and makes  the  moduli space into an integrable system.  In what follows we shall restrict our attention to $GL(n,\C)$-Higgs bundles, which are those Higgs bundles introduced in the first paragraph of these notes, and whose Hitchin fibration is depicted in Figure \ref{fibration}.

 The generic fibre of the Hitchin fibration --- appearing in violet in Figure \ref{fibration} --- is an abelian variety, leading to what is refer to as the \textit{abelianization} of the moduli space of Higgs bundles, and which can be seen geometrically by considering eigenvalues and eigenspaces of the Higgs field. Indeed, a Higgs bundle  $(E,\Phi)$ defines a ramified cover of the Riemann surface given by its eigenvalues and obtained through its characteristic equation:
 \begin{eqnarray}S=\{{\rm det}(\Phi-\eta)=0\} \subset {\rm Tot}K.\label{spectral}\end{eqnarray}
 %
  This cover allows one to construct the  \textit{spectral data} associated to $(E,\Phi)$, which provides a geometric description of the moduli space of Higgs bundles, and is given by
\begin{itemize}
\item the spectral curve $\pi:S\rightarrow \Sigma$, defining a point in the Hitchin base, since the coefficients of $\{{\rm det}(\Phi-\eta)=0\}$ give  a basis of invariant polynomials for the Lie algebra,  and
\item a line bundle  on $S$, defining a point in the Hitchin fibre and obtained as the eigenspace of $\Phi$. 
\end{itemize}
For classical Higgs bundles, the smooth fibres are ${\rm Jac}(S)$, and for $\eta$ the tautological section of $\pi^*K$, one   recovers   $(E,\Phi)$ up to isomorphism from the curve $(S, L\in {\rm Jac}(S))$ by taking $E=\pi_*L$ and $\Phi=\pi_* \eta$.

 When considering $\GC$-Higgs bundles, one has to require  appropriate conditions on the spectral curve and the line bundle reflecting the nature of   $\GC$. This approach originates in the work of Hitchin and of Beauville, Narasimhan and Ramanan (see \cite{PCMI} for references), and we shall describe here an example to illustrate the setting. 
Consider $SL(n,\C)$-Higgs bundles, for which the coefficient in  \eqref{spectral} corresponding to ${\rm Tr}(\Phi)=0$, and the generic fibres of the Hitchin fibration are isomorphic to Prym varieties ${\rm Prym}(S,\Sigma)$.
\smallbreak

\noindent \textbf{Example 3.}~For rank two Higgs bundles, we return to the example in the previous page in which  the Hitchin fibration is over $ H^0(\Sigma, K^2)$, and   the Hitchin map is 
$h:(E,\Phi)\mapsto {\rm det}(\Phi).$ Then,   the family $(E,\Phi_a)$ gives a section of the Hitchin fibration: a smooth map from the Hitchin base to the fibres,  known as \textit{the Hitchin section}. Moreover, this comprises a whole component of the moduli space of real $SL(2,\R)$-Higgs bundles seen inside $\MGC$, and provides an example of a {\it Hitchin component}, also called a {\it Teichmuller component} of real Higgs bundles.\smallskip

In the early 90's Hitchin showed that these components for higher rank split groups $G$,  so-called  {\it Hitchin components}, are homeomorphic to a vector space of dimension $dim(G)(2g - 2)$ and conjectured that they should  parametrize geometric structures --- recall that the Teichm\"uller space $\mathcal{T}(S)$ of the underlying surface $S$ of   $\Sigma$ is the space of marked conformal classes of Riemannian metrics on $S$. These spaces presented the first family of {\it higher Teichm\"uller spaces} within the Betti moduli space of reductive surface group representations $Rep^+(\pi_1(\Sigma),G)$, which leads us to    applications of Higgs bundles within {\it Higher Teichmuller Theory}. 

 \newpage
 \section*{Higher Teichm\"uller theory}

 The Hitchin component of $G$-Higgs bundles, or equivalently of surface group representations,   can be
defined  as the connected component of the Betti moduli space $Rep^+(\pi_1(\Sigma),G)$ containing Fuchsian representations in $G$, which are representations obtained by composing a discrete and faithful representation $\rho:\pi_1(\Sigma) \rightarrow PSL(2, \R)$ with the unique (up to conjugation) irreducible representation $PSL(2, \R)\rightarrow G$. Moreover, as mentioned in the previous section, these representations, called {\it Hitchin representations},  are considered the first example of higher Teichm\"uller space for surfaces: a subset of the set of representations of discrete
groups into Lie groups of higher ranks consisting entirely of discrete and faithful elements.  In order to give the above geometric description of Hitchin representations, 
  Labourie introduced the notion of {\it Anosov representations}, 
    which can be thought of  as a generalization of  convex-co-compact representations to Lie
groups $G$ of higher real rank\footnote{For example, for representations in $SL(2,\C)$, these are quasi-Fuchsian representations.}. 

As beautifully described in Wienhard's 2018 ICM Proceedings article   \cite{Anna}, higher Teichm\"uller theory recently  emerged as a new field in mathematics, closely related to Higgs bundles (see also \cite{Fanny,BrianSigma} for further references).   
 There are  two  types of higher Teichm\"uller spaces,  giving   the only known examples of components which consist entirely of Anosov representations  for surfaces:
 \begin{itemize}
 \item[{\bf (I)}] the space of Hitchin representations into a real split simple
Lie group $G$; and
 \item[{\bf (II)}] the space of maximal representations  into a Hermitian Lie group $G$.
 \end{itemize}
 Recall that a representation $\rho : \pi_1(\Sigma) \rightarrow G$ is maximal if it maximizes the {\it Toledo invariant} $T(\rho)$,  a topological invariant defined for any simple  Lie group  $G$ of Hermitian type as
  \begin{eqnarray}
 \frac{1}{2\pi}\int_{\Sigma}f^*\omega
 \end{eqnarray}
 for $\omega$ the 2-form on the Riemannian symmetric space $X$ of $G$ induced through the imaginary part of the $G$-invariant
Hermitian form on $X$, and for $f:\tilde \Sigma \rightarrow X$ the {\it developing map},  any $\rho$-equivariant smooth map.

%

\smallskip
\noindent \textbf{Example 4.} The Toledo invariant can   be expressed in terms of Higgs bundles. For example for $SL(2,\R)$-Higgs bundles $(L\oplus L^*,\Phi)$, the Toledo invariant is  $2 \deg(L)$ and satisfies $0\leq|2 \deg(L)|\leq 4g-4$. Hence, the family  $(E,\Phi_a)$ from Example 1 is maximal. 
\smallskip

 The existence of   spaces other than those in {\bf (I)} and {\bf (II)} with similar properties to Teichm\"uller space  is currently a topic of significant investigation --- examples are the spaces of $\theta$-positive representations conjectured by Guichard-Wienhard to be potential candidates, and recently shown to exist in \cite{BrianInv}. 

Whilst Anosov representations give a  clear link between discrete and faithful representations and geometric structures, there is no known Higgs bundle characterization of Anosov representations, and very little is known about which explicit geometric structures correspond to these spaces. For instance, work of Choi and Goldman shows that the holonomy representations of convex projective structures are the Hitchin representations when $G=PSL(3,\R)$.  

This brings us to one of the fundamental problems in modern geometry: the classification of geometric structures admitted by a manifold $M$.
Recall that a model geometry is a pair $(G, X)$ where $X$ is a manifold ({\it model space})
and $G$ is a Lie group acting transitively on $X$ ({\it group of symmetries}). Then, a  $(G,X)$-structure on  a manifold $M$  is a maximal atlas of coordinate charts on $M$ with values in $X$ such that the transition maps are given by elements of $G$.   When describing the geometric structures arising through Anosov representations, 
 Higgs bundles have become a key tool through their appearance in relation to higher Teichm\"uller spaces. An example of this is how the Hitchin system was fundamental when showing that maximal representations to $PO(2,q)$ give rise to $(G, X)$-manifolds for at least two choices of $X$: when $X$ is the space   of  null geodesics (photons) in a particular Einstein manifold and when $X=\mathbb{P}(\mathbb{R}^{2+q})$  (e.g. see  \cite{BrianSigma}). 
 For an excellent review of the theory of geometric structures, the reader should refer to Kassel's 2018 ICM Proceedings article \cite{Fanny}.

\newpage
\section*{Harmonic metrics}
 Equivariant harmonic maps play an important role in the non-abelian Hodge correspondence mentioned before (and beautifully reviewed in \cite{whatis}), and thus we shall dedicate this section to look into some of the advances made in this direction.  
In our setting, from the work of Corlette and Donaldson,   any reductive representation $\rho: \pi_1(\Sigma)\rightarrow G_\C$ has associated a  $\rho$-equivariant harmonic map $f$ from the universal cover $\tilde{\Sigma}$ of $\Sigma$ to the corresponding symmetric space of  $G_\C$, which in turn defines a Higgs bundle $(E,\Phi)$. Recall that a map $f:\tilde{\Sigma}\rightarrow M$ is called $\rho$-equivariant if $f(\gamma \cdot x)= \rho(\gamma)\cdot f(x)$ for all $x\in \tilde \Sigma$ and $\gamma \in \pi_1(\Sigma)$. Moreover, such a map induces the $\rho$-equivariant map $df$, leading to the \textit{energy density} defined as
 \begin{eqnarray} e(f)=\frac{1}{2}<df,df>: \tilde \Sigma \rightarrow \R,\label{density}
 \end{eqnarray}
 which is also $\rho$-equivariant and descends to $\Sigma$. Then, the \textit{energy} of $f$ is defined as 
 \begin{eqnarray} E(f)=\int _{M}e(f)d {\rm Vol},  \label{energy}
 \end{eqnarray}
which is finite since $\Sigma$ is compact.   The map $f$ is said to be \textit{harmonic} if it is a critical point of the energy functional $E(f)$ in  \eqref{energy}.

Conversely, through the work of Hitchin 
and Simpson, a polystable Higgs bundle admits a  hermitian metric on the bundle such that the associated Chern connection $A$ solves the Hitchin equations  \eqref{equation1}-\eqref{equation2}, and such a metric is called harmonic. Moreover, the harmonic metric induces an irreducible representation $\rho: \pi_1(\Sigma)\rightarrow G_\C$ and a $\rho$-equivariant harmonic map into the corresponding symmetric space, and these  two
directions together give the celebrated non-abelian Hodge correspondence. 
 
 Understanding the  geometric and analytic properties of the harmonic maps arising from Hitchin's equations  \eqref{equation1}-\eqref{equation2} is of significant importance. For instance, one may ask how do those metrics behave at the boundaries of the moduli space, or how do the energy densities of the corresponding harmonic maps at different points of the Hitchin fibration determine each other (the reader may be interested in the reviews \cite{qiongling} and \cite{Laura}, and references therein). \newpage
 
From Hitchin's work,  the moduli space of Higgs bundles has a natural $\C^*$-action
$\lambda \cdot (E,\Phi) = (E,\lambda \Phi),$
whose fixed point sets allow one to study different aspects of the topology and the geometry of the space, as done in \cite{N2} (see  also  \cite{SteveSigma,BrianSigma}).
Moreover, as shown by Simpson, the  the fixed points by this action are complex variations of Hodge structure (VHS). Recall that a VHS is a  $\mathbb{C}^{\infty}$ vector bundle $V$ with decomposition $V = \bigoplus_{p+q=w}V^{p,q}$, and a  flat connection, satisfying the axioms of Griffiths transversality and existence of a polarization. 
 
\begin{figure}[h]
\centering
\includegraphics[width=0.35\textwidth]{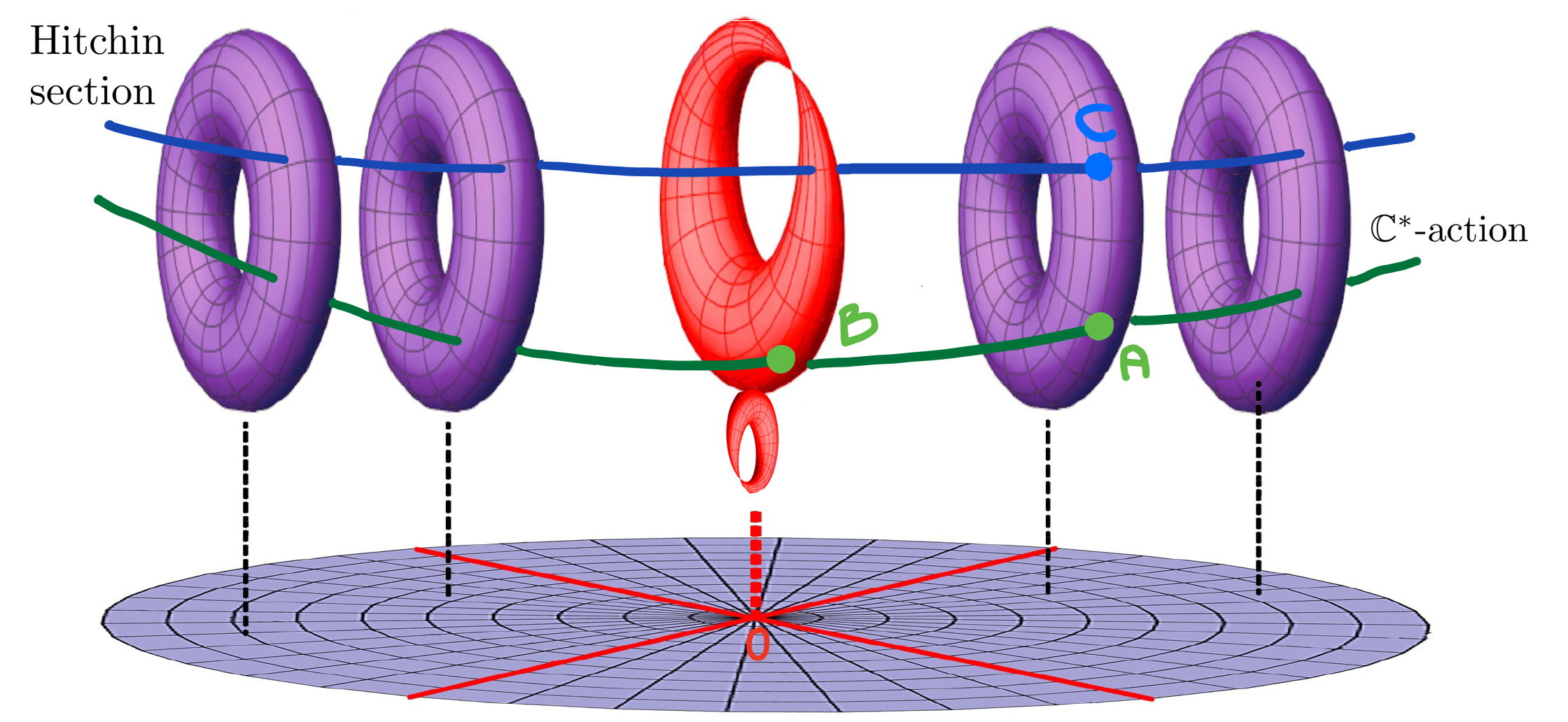}\caption{The nilpotent cone   in red over the 0, and the points $A,B$ and $C$, lying over the fixed point set of the $\C^*$ action and of the Hitchin section respectively.}\label{Li} \end{figure}

From the above, one may ask how the energy density of harmonic maps
changes along the $\C^*$-flow on the moduli space of Higgs bundles. Whilst this remains a challenging open question in the area, a better understanding might come from the following conjectural picture of Dai-Li described  in Figure \ref{Li}, and through which the harmonic map of a fixed point set of the $\C^*$ action on $\MGC$ gives rise to two other related harmonic maps. 

Here, given a point $A$ within the Hitchin fibration,  one can immediately determine
the point $B$ to be the limit of the $\C^*$-flow $\lambda \cdot A$ as $\lambda \rightarrow 0$ in the nilpotent cone, and the point $C$
to be the intersection point of the Hitchin fiber containing $A$ and the Hitchin section. Then Dai-Li's conjecture states that the energy densities defined as in Eq.\eqref{density} of the corresponding harmonic maps
$f_A, f_B, f_C : \tilde{\Sigma} \rightarrow \mathbb{N}$ satisfy
\begin{eqnarray}e(f_B) < e(f_A) < e(f_C).\label{li}\end{eqnarray}
As evidence for the above conjecture, one can consider the integral version (through Eq.\eqref{energy}), for which   Hitchin showed  that $E(f_B) < E(f_A)$, but where the other corresponding inequality in  \eqref{li} remains open. 
\newpage
\section*{Limiting structures}
The study of $\rho$-equivariant harmonic metrics and higher Teichm\"uller theory through Higgs bundles has received much attention in recent years and brings us to one of the most important conjectures in the area. This conjecture, due to Labourie,   states   that for each Hitchin representation $\rho$ there is a unique conformal structure $X_\rho$ on the underlying surface $S$ in which the $\rho$-equivariant harmonic metric is a minimal immersion.
In particular, Labourie showed that for all Anosov representations such a conformal structure exists, but   the difficultly lies in proving uniqueness ---  the conjecture  has been established only for   Lie groups of rank two (\cite{MR3583351,BrianSigma}).  To understand this problem, it becomes fundamental to study the deformation of  conformal structures on surfaces and the corresponding
 harmonic metric. 

Some of these deformations can be see through the hyperk\"ahler structure of the moduli space, though which it has a $\mathbb{C}\mathbb{P}^1$-worth of complex structures labelled by a parameter $\xi\in \C^\times$. Indeed, we can think of a hyperk\"ahler manifolds as a manifold whose tangent space admits an action of three complex structures $I$, $J$ and $K$ satisfying the quaternionic equations and compatible with a single metric. In our case, $I$ arises from the complex structure $I$ on the Riemann surface $\Sigma$, while   $J$ is from the complex structure on the group $\GC$.   In this setting, one has the following spaces:
 \begin{itemize}
 \item   $\xi=0$ gives the moduli space of Higgs bundles, 
 \item    $\xi\in\C^\times$ gives  the moduli space of flat connections\begin{eqnarray}
\nabla_\xi=\xi^{-1}\Phi+\bar\partial^h_A+ \partial^h_A+\xi\Phi^{^*h};\label{Laura1}
\end{eqnarray} 
 \item finally,  taking $\xi=\infty$ gives the moduli space of so-called ``anti-Higgs bundles''.
 \end{itemize}  

The hyperk\"ahler metric on Hitchin moduli space is expected to be of type ``{\it quasi-ALG}'' which is some expected generalization of ALG. Indeed, a far reaching open question is the understanding of the behaviour of  the metrics  at the boundaries of the space. For instance taking the limit  of Hitchin's solutions along a ray in the Hitchin base \[\lim_{t\rightarrow \infty}(\bar\partial_A,t\Phi,h_t).\]   

Almost a decade ago Gaiotto-Moore-Neitzke gave a conjectural description of the hyperk\"ahler metric on $\MGC$, which surprisingly suggests that   much of the asymptotic geometry of the moduli space can be derived from the abelian spectral data described before. Recent progress has been made by Mazzeo-Swoboda-Weiss-Witt, Dumas-Neitzke and Fredrickson but the global picture remains open (for a survey of the area, see \cite{Laura}).

Finally, one further type of limiting structure we would like to mention is that of opers, appearing as certain limits of Higgs bundles in the Hitchin components. To see this,  note that for a solution of   \eqref{equation1}-\eqref{equation2}  in the  $SL(n, \C)$-Hitchin section, one can add a real parameter $R>0$ to \eqref{Laura1} to obtain a natural family of  connections with $SL(n,\R)$ monodromy
\begin{eqnarray}
\nabla(\xi,R):=\xi^{-1}R\Phi+\bar\partial_A+ \partial^h_A+\xi R\Phi^{^*h}.
\end{eqnarray}

 Some years ago Gaiotto  conjectured that the space of opers should be obtained as the $\hbar$-conformal limit of the Hitchin section: taking   $R \rightarrow 0$ and  $\xi\rightarrow 0$ simultaneously while
holding the ratio $\hbar= \xi/R$ fixed. 

 The conjecture was recently established for general simple Lie groups by Dumistrescu-Fredrickson-Kydonakis-Mazzeo-Mulase-Neitzke, who also conjectured that that this oper is the {\it quantum curve} in the sense of Dumitrescu-Mulase, a quantization of the spectral curve  $S$ of the corresponding Higgs bundle by {\it topological recursion} (see references and details in \cite{Olivia}). 
More recently,
Collier-Wentworth showed that  the above conformal limit exists in much more generality and gives a correspondence between (Lagrangian) strata for every stable VHS --- and not only the Hitchin components. Specifically, they constructed a generalization of the Hitchin section by considering stable manifolds $\mathcal{W}^0(E_0,\Phi_0)$ arising from each VHS $(E_0,\Phi_0)$ given by
\begin{eqnarray}
\{(E,\phi)\in \MGC ~|~ \lim_{t\rightarrow 0} t\cdot (E,\Phi)=(E_0,\Phi_0)\}.
\end{eqnarray} 
The equivalent correspondence is then obtained by identifying a subcomplex of the deformation complex at a VHS with a corresponding slice in the space of Higgs bundles which parametrizes $\mathcal{W}^0(E_0,\Phi_0)$.

\newpage
\section*{Correspondences}

%

The appearance of Higgs bundles as parameter spaces for geometric structures is an example of the study of correspondences between solutions to Hitchin's equations  \eqref{equation1}-\eqref{equation2} and different mathematical objects. 
In what follows we shall restrict our attention to a few correspondences   between Higgs bundles and two classes of  mathematical objects: quiver varieties and hyperpolygons (e.g.~see references in  \cite{VicSigma,SteveSigma}). 

\begin{figure}[h!]
\centering
\includegraphics[width=0.18\textwidth]{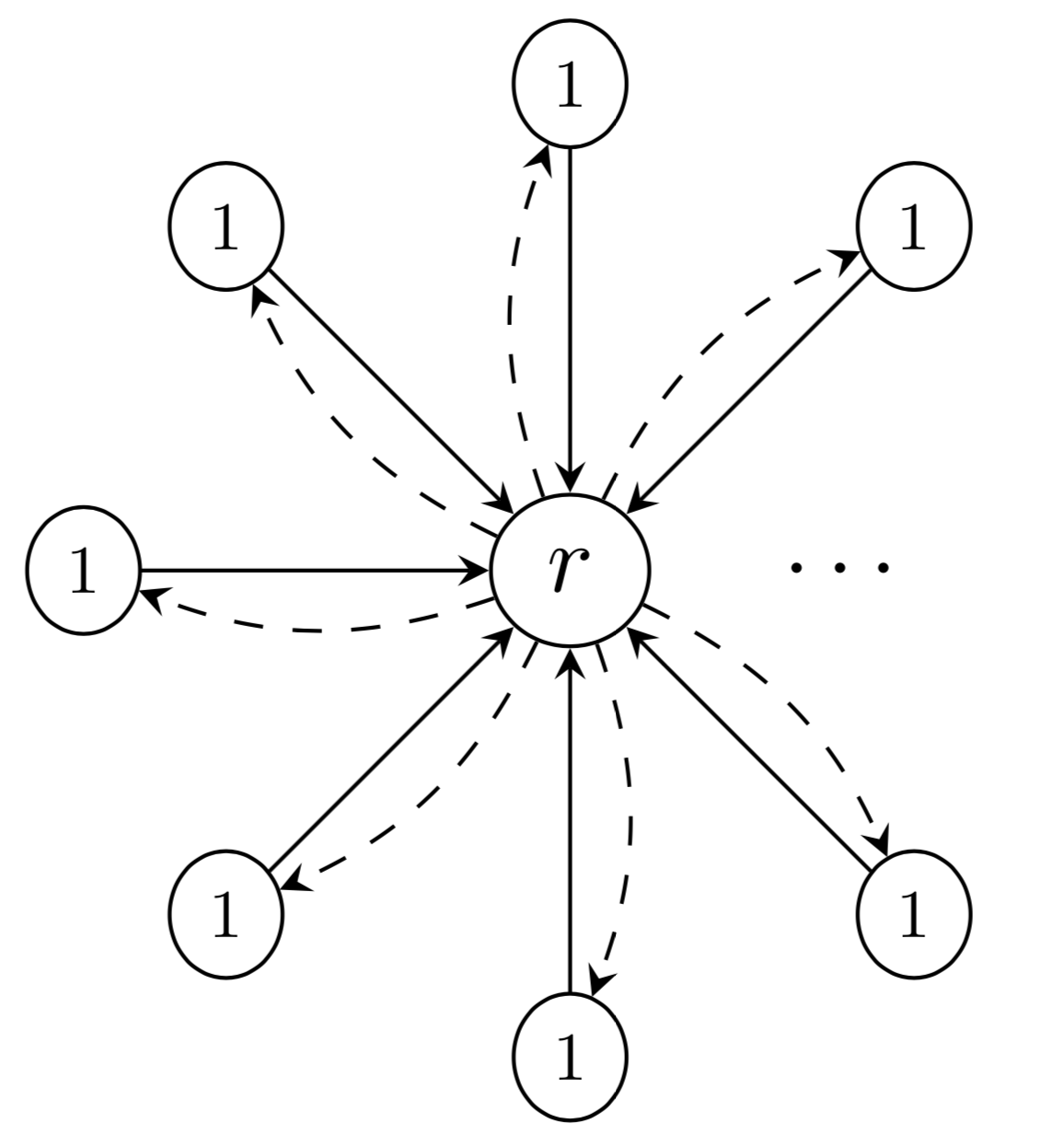}  \caption{A star-shaped quiver. 
}\label{hyper}\end{figure}

Recall that 
a quiver $Q=(V,A,h,t)$ is an oriented graph, consisting of a finite vertex set $V$, a finite arrow set $A$, and head and tail maps $h,t:A\rightarrow V$.
A Nakajima representation of a quiver $Q$ can   be written as families  
$W:=((W_v), \phi_a, \psi_a)$ for $a\in A$ and  $v\in V$, where $W_v$ is a finite dimensional vector space. The map
   $\phi_a: W_{t(a)}\rightarrow W_{h(a)}$ is a linear map for all $a\in A$, and $\psi_a$  is in the
cotangent space to
$Hom(W_{t(a)},W_{h(a)})$
at $\phi_a$.  In particular, a \emph{hyperpolygon} is a representation of the \emph{star-shaped} quiver, an example of which  appears in Figure \ref{hyper}.

%
%
%
%
%
%
%
%
%

 For the star-shaped quiver in Figure \ref{hyper}, for which the   dimensions of $W_v$ are indicated in each vertex, the vector space $T^*{\rm Rep}(Q)$ of representations of $Q$  is
\begin{eqnarray}
 T^* {\rm Rep}(Q) &=& T^*\left(\bigoplus_{i=1}^n {\rm Hom}(\C, \C^r)\right)\nonumber\\
&\cong&   T^* \left({\rm Hom}(\C^n, \C^r)\right) =T^*( \C^{n\times r}).\nonumber
\end{eqnarray}

Konno showed that hyperpolygon spaces are   hyperk\"ahler analogs of polygon spaces. Moreover, through the work of Fisher-Rayan, the space of hyperpolygons as in Figure \ref{hyper} may be identified with a moduli space of certain rank $r$ parabolic Higgs bundles on $\mathbb{P}^1$. 
In this setting, one has to puncture $\mathbb{P}^1$ along a positive divisor $D$ and then regard the Higgs field as being valued in $\mathcal{O}(q) = K\otimes \mathcal{O}(D)$, with poles along $D$ and satisfying  certain conditions on its residues  at the poles. This takes us to a  generalization of Higgs bundles on higher genus surfaces  obtained by allowing the Higgs field to have poles, leading to the moduli spaces of tame or parabolic Higgs bundles (for logarithmic singularities) initiated by Simpson \cite{simpson92}, or of wild Higgs bundles (for higher order poles) initiated by Boalch --- see references in  \cite{SIGMA} to learn more about these other settings.
Understanding the more general appearance of parabolic (and wild) Higgs bundles on higher genus Riemann surfaces in terms of hyperpolygons remains an open question.

In a different direction, given a fixed Riemann surface $\Sigma$ and a homomorphism between two Lie
groups $\Psi : G_\C \rightarrow G_\C',$ there is a naturally induced map between representation spaces modulo conjugation
\[\Psi : Rep(\Sigma, G_\C) \rightarrow Rep(\Sigma, G_\C'),\]
It follows from the non-abelian Hodge correspondence  that there must be a corresponding induced map between the Higgs bundle moduli spaces,  but this 
 does not transfer readily to the induced map on Higgs bundles, in particular since the image might be over the singular locus of the base. 
  When the maps arise through isogenies, together with Bradlow and Branco, the author obtained a description of the map for spectral data in terms of fibre products of spectral curves \cite{PCMI}, but of much interest is the understanding of other maps arising in this manner.

\begin{figure}[h]
\centering
\includegraphics[width=0.34\textwidth]{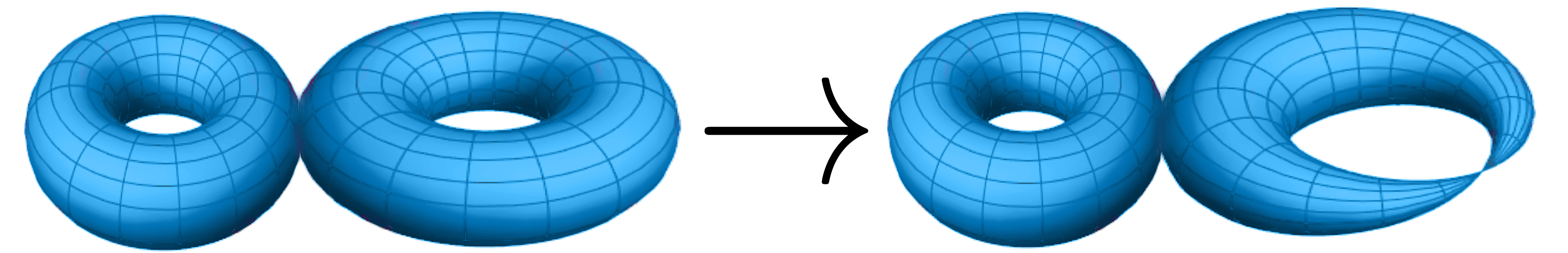} \caption{A degeneration of the Riemann surface.}\label{deg}\end{figure}

Finally,   when considering compactifications of the moduli space, one may   ask how do the moduli spaces transform when the base Riemann surface $\Sigma$ changes (for instance, when  degenerating the surface $\Sigma$ as in Figure \ref{deg}), a question closely related to the relation between Higgs bundles and singular geometry, which we shan't touch upon here --- see \cite{SIGMA} for a survey and open problems in this direction. 
 
\newpage

\section*{Mirror symmetry and branes}

One of the most interesting correspondences of Hitchin systems arises through mirror symmetry. For $^L G_\C$   the Langlands dual group of $G_\C$, there is   an identification  of the Hitchin basis $\mathcal{A}_{G_\C} \simeq \mathcal{A}_{^L G_\C}$. Then, through the famous SYZ conjecture, mirror symmetry should manifest as a duality between the   spaces of Higgs bundles for Langlands dual groups   fibred over the same   base via the Hitchin fibration.  The two moduli spaces $\mathcal{M}_{G_\C}$ and $\mathcal{M}_{^L G_\C}$ are then torus fibrations over a common base and as first observed by Hausel-Thaddeus for $SL(n,\C)$ and $PGL(n,\C)$, and shown by Donagi and Pantev for general pairs of Langlands dual reductive groups, their non-singular fibres are dual abelian varieties. Kapustin and Witten gave a physical interpretation of this in terms of S-duality, using it as the basis for their approach to the geometric Langlands program.

   The appearance of Higgs bundles (and flat connections) within string theory and the geometric Langlands program has led researchers to study the {\it derived category of coherent sheaves} and the {\it Fukaya category} of these moduli spaces. Therefore, it has become fundamental to understand Lagrangian submanifolds of the moduli space of Higgs bundles supporting holomorphic sheaves ($A$-branes), and their dual objects ($B$-branes). 
      By considering the support of branes, we shall refer to a submanifold of a  hyperk\"ahler manifold as being of type $A$ or $B$ with respect to each of the complex   structures $(I,J,K)$.
       Hence one may study branes of the four possible types: $(B,B,B), (B,A,A), (A,B,A)$ and $(A,A,B)$, whose dual partner is predicted by Kontsevich's homological mirror symmetry to be:
      \begin{eqnarray}
       (B,A,A) \longleftrightarrow (B,B,B) \label{uno}\\ (A,A,B) \longleftrightarrow (A,A,B)\\
              (A,B,A) \longleftrightarrow (A,B,A)        \end{eqnarray}

In views of the  SYZ conjecture,   it is crucial to obtain the duality between  branes within
the Hitchin fibration, and in particular between those completely contained within the irregular fibres, and   this has remained a very fruitful direction of research for decades. 
    In 2006  Gukov, Kapustin and Witten  introduced the first studies of branes of Higgs bundles in relation to the  Geometric Langlands Program and electromagnetic duality where the $(B,A,A)$-branes of real $G$-Higgs bundles were considered.  These branes, which correspond to surface group representations into the real Lie group $G$, may intersect the regular fibres of the Hitchin fibration in very different ways (see \cite{OBW,PCMI} for references):
   
   \begin{itemize}
   \item {\it Abelianization --- zero-dimensional intersection.} When $G$ is a split real form, the author showed that the $(B,A,A)$-brane intersects the Hitchin fibration in torsion two points. 
   
     \item {\it Abelianization --- positive dimensional intersection.} Moreover, we can also show that for other groups such as $SU(n,n)$, the intersection has positive dimension but may still be described solely  in terms of abelian data.  
   \item {\it Cayley/Langlands type correspondences.} Surprisingly, many     spaces of Higgs bundles   corresponding to non-abelian real gauge theories  do admit abelian parametrizations via  auxiliary spectral curves,  as     shown with Baraglia    through Cayley/Langlands type correspondences for the groups $G=SO(p+q,p)$ and $G=Sp(2p+2q,2p)$.

   \item  {\it Nonabelianization. }Together  with Hitchin we initiated the study of branes which don't intersect the regular locus, through the  {\it nonabelianization of Higgs bundles}, which characterized   the branes for $G=SL(n,\mathbb{H}),~SO(n,\mathbb{H})$ and $Sp(n,n)$ in terms of non-abelian data given by   spaces of rank 2 vector bundles on the spectral curves.

   \end{itemize}
   
   Moreover, it has been conjectured (Baraglia-Schaposnik) that the Langlands dual in \eqref{uno} to the above $(B,A,A)$-branes  should correspond to the  $(B,B,B)$-branes of Higgs bundles with structure group the {\it Nadler group} \cite{OBW}. More generally, branes of Higgs bundles have shown to be
 notoriously difficult to compute in practice, and very few broad classes of examples are known --- e.g. see \cite{PCMI} for a partial list of examples. In the next section  we shall describe a family of branes obtained by the author and Baraglia
 by imposing symmetries to the solutions of \eqref{equation1}-\eqref{equation2} --- see \cite{OBW} and references therein.

 \newpage

\section*{Higgs bundles and 3-manifolds}
By considering actions on the Riemann surface $\Sigma$ and on the Lie group $\GC$, one can induce actions on the moduli space of Higgs bundles and on the Hitchin fibration, and study their fixed point sets. Indeed, together with Baraglia we defined  the following:
\begin{itemize} 
\item  Through the Cartan involution $\theta$ of a real form $G$ of $G_\C$ one obtains  
$
i_1(\bar \partial_A, \Phi)=(\theta(\bar \partial_A),-\theta( \Phi)).
$
\item  A real structure  $f : \Sigma \to \Sigma$ on $\Sigma$  induces 
$
i_2(\bar \partial_A, \Phi) = (f^*(\rho(\bar \partial_A)), -f^*( \rho(\Phi) )).
\label{i2}
$
\item  Lastly, by looking at $i_3 = i_1 \circ i_2$, one obtains   
$
i_3(\bar \partial_A, \Phi)=(f^* \sigma(\bar \partial_A),f^*\sigma( \Phi)).
$
\end{itemize}

\begin{figure}[h]
\centering
\includegraphics[width=0.45\textwidth]{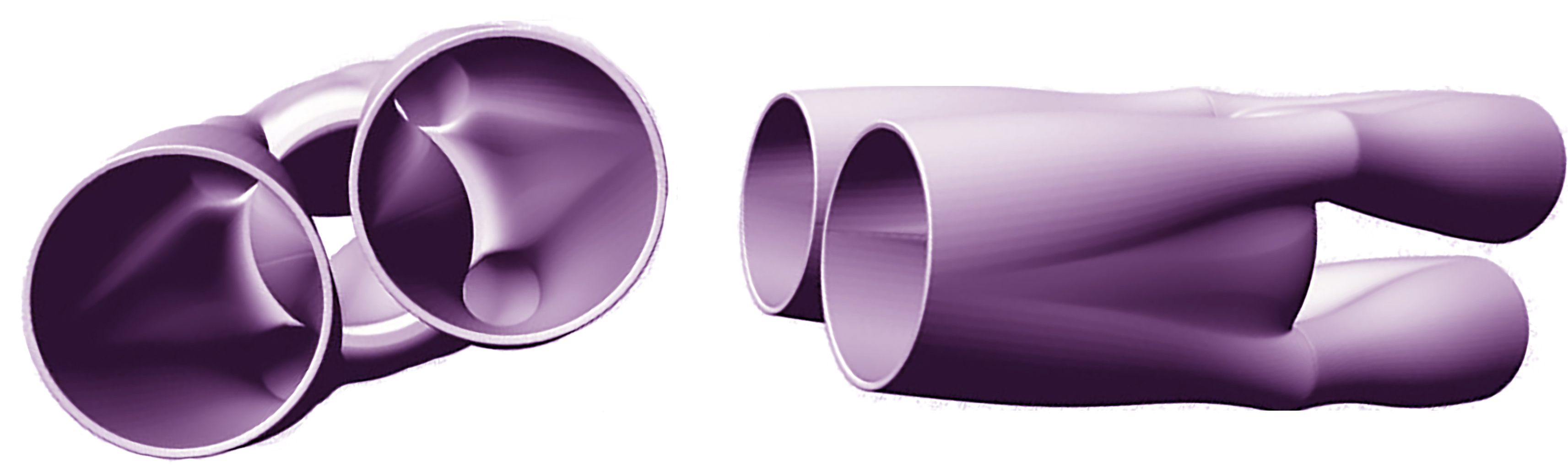} \caption{A real slice fixed by $i_1$ of the moduli space of $SL(2,\C)$-Higgs bundles, from two different angles, obtained through Hausel's 3d prints of slices of $\MGC$.} \label{realslice}\end{figure}

The fixed point sets of $i_1,i_2,i_3$ are branes of type $(B,A,A),(A,B,A)$ and $(A,A,B)$ respectively. The topological invariants can be described using $KO$, $KR$ and equivariant $K$-theory \cite{OBW}. In particular,  the fixed points of $i_1$  give the $(B,A,A)$-brane of $G$-Higgs bundles mentioned in the previous section,  an example of which appears in Figure \ref{realslice}, and which one can study via the monodromy action on the Hitchin fibration (e.g. see \cite{PCMI}).

  The fixed points of the involution $i_2$ are representations of the orbifold fundamental group of certain 3-manifold obtained through $f$ bounding $\Sigma$. Recall that a real structure (or  anti-conformal maps) on a compact connected Riemann surface $\Sigma$ is an 
 anti-holomorphic involution $f : \Sigma \to \Sigma$.  
 
The classification of real structures on compact Riemann surfaces  is a classical result of Klein, who showed that all such involutions on $\Sigma$ may be characterised by two integer invariants $(n,a)$: the number $n$ of   disjoint union of   copies of the circle embedded in $\Sigma$ fixed by $f$; and $a\in \mathbb{Z}_2$  determining whether the complement of the fixed point set has one ($a=1$) or two ($a=0$) components, e.g. see 
Figure \ref{hola}.

 \begin{figure}[h]
\centering
\includegraphics[width=0.45\textwidth]{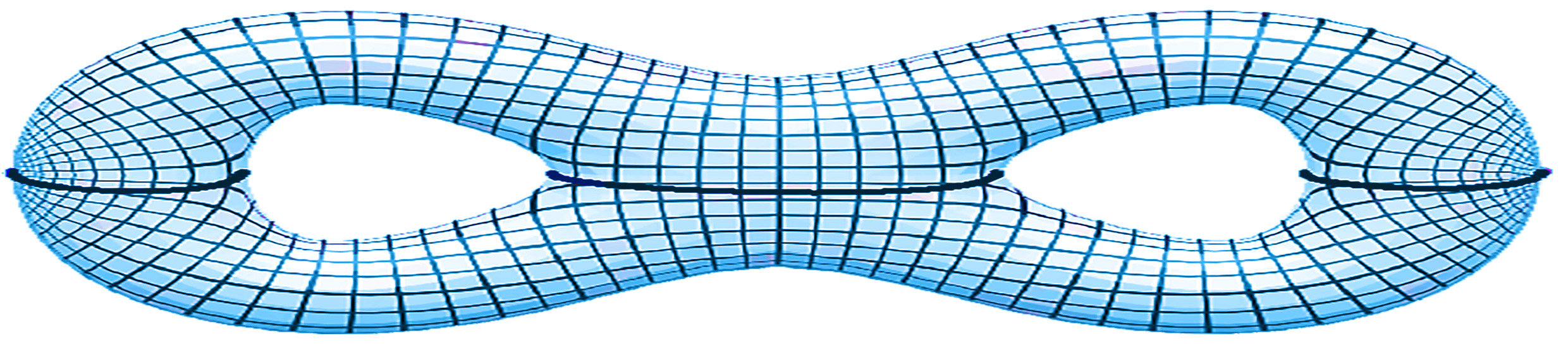} \caption{A genus 2 Riemann surface and its fixed point sets under an anti-holomorphic involution with invariants $(n,a)=(3,0)$.}\label{hola}\end{figure} 

A real structure $f$ on the Riemann surface $\Sigma$ induces  involutions on the moduli space of representations $\pi_1(\Sigma) \to \GC$, of flat connections and  of $\GC$-Higgs bundles on $\Sigma$, and the fixed points sets define $(A,B,A)$-branes of Higgs bundles. These branes can be shown to be real integrable systems, given by (possibly singular) Lagrangian fibrations. 

From a representation theoretic point of view, one may ask which interesting representations these branes correspond to, a question closely related to the understanding of \textit{which representations of $\pi_1(\Sigma)$ extend to $\pi_1(M)$ for $M$ a 3-manifold whose boundary is $\Sigma$}. Whilst this question in its full generality remains an important open problem, we can consider some particular cases in which the answer becomes clear from the perspective of Higgs bundles. For this, as seen in \cite{OBW} and references therein, we considered  the 3-manifolds  \begin{eqnarray}M= \frac{\Sigma \times [-1,1]}{ (x,t)\mapsto (f(x),-t)},\label{mani}\end{eqnarray} for which $\partial M = \Sigma$   (e.g. a handle body).
In this setting, together with Baraglia, we were able to show that a connection solving the Hitchin equations \eqref{equation1}-\eqref{equation2} on $\Sigma$ extends over $M$ given in \eqref{mani} as a flat connection iff the  Higgs bundle $(E,\Phi)$ is fixed by $i_2$ and the class $[E] \in \tilde{K}^0_{\mathbb{Z}_2}(\Sigma)$ in reduced equivariant $K$-theory is trivial. This is, the Higgs bundles which will extend are only those whose vector bundle is preserved by the list of the involution $i_2$, and for which the action of $i_2$ over the fibres of $E$  is trivial when restricted to each fixed circle (as those circles in Figure \ref{hola}).


\newpage

 \section*{Global topology}

The computation of topological invariants of Higgs bundle moduli spaces has received vast attention from researchers who have tackled this problem with a diverse set of mathematical tools --- see \cite{Hausel2011GlobalTO} and \cite{SteveSigma} and references therein. One of the main questions considered for Higgs bundles and their generalizations is what the Poincar\'e polynomial $P(\MGC,t)$ of the space is. 
A useful fact is that the total space of the Hitchin fibration deformation retracts onto the nilpotent cone $h^{-1}(0)$ via the gradient flow of the moment map of the $C^*$-action introduced in the harmonic metrics section. Hence,  the cohomology ring localises to the fixed-point locus inside $h^{-1}(0)$: as first seen in the work of Hitchin, the Poincar\'e series  that generates the Betti numbers of the rational cohomology $H^\bullet(\MGC, \mathbb{Q})$ are a weighted sum of the Poincar\'e series  of the connected components of the fixed-point locus, which is essentially Morse theory. 

 \smallbreak
 
 \noindent{\bf Example 5.} As shown by Hitchin, for the $SL(2,\C)$-Higgs bundles  in {\it Example 1}, when the genus of $\Sigma$ is $g = 2$, the Poincar\'e series is
 \begin{eqnarray}
 1+t^2+4t^3+2t^4+34t^5+2t^6.
 \end{eqnarray}
 \smallbreak
 
  Using Morse theory, it has only been possible to compute Poincar\'e polynomials for low rank groups, and extending this to higher rank has been a challenging open problem for some time.  More recently,  interesting alternative techniques have been used to access the higher-rank Poincar\'e polynomials by Mozgovoy, Schiffmann, Mellit, and others.
 One may further ask about the  structure of the ring $H^\bullet(\MGC, \mathbb{Q})$ itself: for instance Heinloth recently proved that the intersection pairing in the middle dimension for the smooth moduli space vanishes in all dimensions for $G_\C=PGL(n,\C)$; and Cliff-Nevins-Shen proved that that the Kirwan map from the cohomology of the moduli stack of $G$-bundles to the moduli stack of semistable $G$-Higgs bundles fails to be surjective.

One of the most important cohomological conjectures in the area is de-Cataldo-Hausel-Migliorini's {\it P=W conjecture}, which gives a correspondence between the  weight filtration and the perverse filtration on the cohomology of $\mathcal{M}_B$ and $\mathcal{M}_{Dol}$, respectively, obtained via  non-abelian Hodge theory. Only certain special cases are known, e.g., for rank 2 Higgs bundles, shown by de-Cataldo-Hausel-Migliorini's (see \cite{Hausel2011GlobalTO}), and for certain moduli spaces of wild Higgs bundles, proven recently by Shen-Zhang and Szabo.

 Inspired by the SYZ conjecture mentioned before,  Hausel-Thaddeus conjectured that
mirror moduli spaces of  Higgs bundles present an agreement of
appropriately defined  Hodge numbers:
\begin{eqnarray}
h^{p,q}(\MGC )=h^{p,q} (\mathcal{M}_{^LG_\C} ).\label{hola}
\end{eqnarray}
Very recently, the first proof of this conjecture was established  for the moduli spaces of  $SL(n,\C)$ and $PGL(n,\C)$-Higgs bundles by Groechenig-Wyss-Ziegler in   \cite{groechenig2017mirror}, where they  established the equality of stringy Hodge numbers using p-adic integration relative to the fibres of the Hitchin fibration, and interpreted canonical gerbes present on these moduli spaces as characters on the Hitchin fibres.

Further combinatorial properties of $\MGC$ can be glimpsed through their twisted version, consisting of Higgs bundles $(E, \Phi)$ on $\Sigma$ with $\Phi: E\to E\otimes\mathcal L$, where $\Sigma$ now has any genus, $L$ is a line bundle with $\deg(L)>\deg(K)$,    but without any punctures or residues being fixed. 
The corresponding moduli spaces     carry a natural $\mathbb{C}^*$-action but  
but are not hyperk\"ahler and there is no immediate relationship to a character
variety.  Hence, there is no obvious reason for the Betti numbers to be invariant with regards to the choice of deg(E), which would normally follow from non-abelian Hodge theory. However, the independence holds in direct calculations of the Betti numbers in low rank, and was recently shown for $GL(n,\C)$ and $SL(n,\C)$-Higgs bundles  by Groechenig-Wyss-Ziegler in   \cite{groechenig2017mirror}. This suggests that some topological properties of Hitchin systems are independent of the hyperk\"ahler geometry (see references in \cite{Hausel2011GlobalTO,SteveSigma} for more details).

Finally, it should be mentioned that an alternative description of the Hitchin fibration can be given through Cameral data, as introduced by Donagi and Gaitsgory, and this perspective presents many advantages, in particular when considering correspondences   arising from mirror symmetry and Langlands duality, as those mentioned in previous sections studied by Donagi-Pantev.
\newpage

\section*{Acknowledgements} The author is partially supported by the
NSF grant DMS-1509693, the NSF CAREER Award DMS-1749013, and by
the Alexander von Humboldt Foundation. This material is  based upon
work supported by the NSF under Grant No. DMS1440140 while the author was in residence at the Mathematical Sciences Research Institute in Berkeley, California, during the Fall 2019 semester. Finally, the author would like to thank  Steve Bradlow, Laura Fredrickson, Tamas Hausel, Nigel Hitchin,  Qiongling Li,  Sara Maloni, Du Pei, James Unwin, Richard Wentworth,  Anna Wienhard, and very especially Steve Rayan for useful comments on a draft of the manuscript.

\bibliography{Schaposnik_Laura_July2019.bib}

\end{document}